\documentclass{article}
\usepackage{amsmath}
\usepackage{latexsym,amssymb}
\usepackage{amsmath}
\usepackage{stmaryrd}
\usepackage{hyperref}
\usepackage{pifont}
\usepackage{relsize}
\usepackage{lipsum}
\usepackage{tikz}
\usepackage{graphicx}
\usepackage{lmodern,graphicx}
\usepackage{lmodern}
\usepackage{scrextend}

\newtheorem{ex}{Example}[section]
\newtheorem{defn}[ex]{Definition}
\newtheorem{prop}[ex]{Proposition}
\newtheorem{thm}[ex]{Theorem}
\newtheorem{cor}[ex]{Corollary}
\newtheorem{lem}[ex]{Lemma}
\newtheorem{rem}[ex]{Remark}
\newenvironment{pf}{{\bf Proof:} }{$\Box$
\mbox{}}

\def \A {{\mathbf A}}
\def \B {{\mathbf B}}
\def \X {{\cal X}}
\def \G {{\cal G}}

\def \xra{\xrightarrow}
\def \r {\vartriangleleft}
\def \l {\vartriangleright}
\def \d {\partial}
\def \k {\kappa}
\def \ll {\llbracket}
\def \rr {\rrbracket}
\def \la {\l_{\scriptsize{\vdash}}}
\def \lb {\l_{\scriptsize{\dashv}}}
\def \ra {\r_{\scriptsize{\vdash}}}
\def \rb {\r_{\scriptsize{\dashv}}}

\def \C {\mathcal{C}}

\input xypic
\input xy
\xyoption{v2}
\xyoption{all}
\xyoption{2cell}
\xyoption{rotate}

\begin{document}

\title{From Simplicial Homotopy to Crossed Module Homotopy in Modified Categories of Interest}

\author{
	Kad\.{i}r Em\.{i}r\thanks{Corresponding author. The authors are thankful to Enver \"Onder Uslu for his comments on the paper.} \\ {\it  \small k.emir@campus.fct.unl.pt} \\ {  \small Centro de Matem\'atica e Aplica\c{c}\~oes,}\\ { \small Faculdade de Ci\^encias e Tecnologia}{ \small (Universidade Nova de Lisboa), Caparica,  Portugal.} \\
	{ \small Also at: Department of Mathematics and Computer Science,}\\ { \small  Eski\c{s}ehir Osmangazi University, Turkey.}
	\\ \quad \\
	Selim \c{C}etin \\ {\it  \small selimc@ogu.edu.tr} \\ {  \small Department of Mathematics and Computer Science,}\\ { \small Eski\c{s}ehir Osmangazi University, Turkey.}}

\date{\today}
\maketitle

\begin{abstract}
We address the (pointed) homotopy of crossed module morphisms in modified categories of interest; which generalizes the groups and various algebraic structures.  We prove that, the homotopy relation gives rise to an equivalence relation; furthermore a groupoid structure, without any restriction on neither domain nor co-domain of the crossed module morphism. Additionally, we consider the particular cases such as associative algebras, Leibniz algebras, Lie algebras and dialgebras of crossed modules of this generalized homotopy definition. Then as the main part of the paper, we prove that the functor from simplicial objects to crossed modules in modified categories of interest preserves the homotopy and also the homotopy equivalence. 
\end{abstract}

AMS 2010 Classification: {55U10 (principal), 
18D05,    
55P10,  
(secondary). 



\section{Introduction}

Categories of interest were introduced to unify definitions and properties
of different algebraic categories and different algebras. The first steps
for this unification were given by P. G. Higgins in \cite{zbMATH03121671} under the name of ``Groups with multiple operators". (for details, see \cite{zbMATH03957392}). Then the generalized notion ``Categories of interest" was introduced by M. Barr and G. Orzech in \cite{Orz}. Categories of groups, Lie algebras,
Leibniz algebras, (associative) commutative algebras, dialgebras and many
others are basic examples of categories of interest. Nevertheless, the cat$%
^{1}$-algebras are not categories of interest. These categories with a
modification in one condition was introduced in \cite{BDCU} and called it
``Modified categories of interest" which will be denoted by MCI hereafter. Cat%
$^{1}$--Lie (Leibniz, associative, commutative) algebras and many others or
crossed modules of algebras are all MCI \cite{zbMATH05057324, zbMATH06168638,  jas, zbMATH06271154, zbMATH04075344, CDL4} but they are not categories of interest. 

\bigskip 

	The categories \ $\mathbf{Cat}^{\mathbf{1}}$-$\mathbf{Ass}, \ \mathbf{Cat}%
	^{\mathbf{1}}$-$\mathbf{Lie}, \ \mathbf{Cat}^{\mathbf{1}}$-$\mathbf{Leibniz}%
	,\ \mathbf{PreCat}^{\mathbf{1}}$-$\mathbf{Ass}, \ $ $\mathbf{PreCat}%
	^{\mathbf{1}}$-$\mathbf{Lie}$ and $\mathbf{PreCat}^{\mathbf{1}}$%
	-$\mathbf{Leibniz}$ are MCI, which are not
	categories of interest. Also the category of commutative von Neumann regular
	rings is isomorphic to the category of commutative rings with a unary
	operation $(\ )^{\ast}$ satisfying two axioms, defined in \cite{jan}, which is a MCI.
 
 \medskip
 
 A crossed module \cite{BHS} $\G=(\d\colon E \to G, \l)$ of groups, is given by a group homomorphism $\d\colon E \to G$, together with an action $\l$ of $G$ on $E$, such that satisfying the following Peiffer-Whitehead relations for all $e,f \in E$ and $g \in G$:
 \begin{align*}
 \textrm{ \bf{First Peiffer-Whitehead Relation (for groups)}: }\, \d(g \l e)& =g \,\d(e) \, g^{-1},&&&&& &&&&&\\  \textrm{\bf{ Second Peiffer-Whitehead Relation (for groups)}: }\, \d(e)  \l f&=e\, f \, e^{-1}. &&&&& &&&&&
 \end{align*}
 
 Crossed modules were introduced for groups by Whitehead \cite{W3,W4} as algebraic models for homotopy 2-types \cite{B1,Lo}. Another result is that, the category of crossed modules are also equal to cat$^1$ groups \cite{Lo}; therefore to the categories of interest in the sense of \cite{CL4, CL2}. However since the category of some cat$^1$ algebras are not category of interest but are MCI, we will work on this modified category in this paper. In MCI, notion of the crossed module notion introduced in \cite{BDCU}. Crossed modules are also appear in the context of simplicial homotopy theory, since they are  equivalent to simplicial objects with Moore complex of length one in (modified) categories of interest \cite{YO1, C1} which can be diagrammed by:
 \begin{align}\label{equality}
 \xymatrix@C=15pt@R=30pt{
 	Simp(\C)_{\leq 1} \ar[rr]^{X_1} & & XMod \ar[dl] \\
 	& Tr_{1}Simp(\C)_{\leq 1} \ar[ul]^{t_1}}
 \end{align}
 
 \medskip
 
 An equally well established result of this equivalence is that the homotopy category of $n$-types is equivalent to the homotopy category of simplicial groups with Moore complex of length $n-1$, also called algebraic models for $n$-types.
 
\medskip

The homotopy relation between (pre)crossed module morphisms $\G \to \G'$ was introduced for groups by J.Faria Martins in \cite{JFM1}, and for commutative algebras in \cite{IJK}. In both of these studies, we see that the homotopy relation between crossed module morphisms $\G \to \G'$ is an equivalence relation in the general case, with no restriction on $\G$ or $\G'$. If we examine this result in the sense of \cite{DS}, this is an unexpected situation indeed, since the homotopy relation of morphisms $\G \to \G'$ gives an equivalence relation when  $\G =(\d\colon E \to G, \l)$ is cofibrant. On the other hand, \cite{CG},  $\G =(\d\colon E \to G, \l)$  is cofibrant if, and only if,  $G$ is  a free group in the well known model category structure, in the sense of \cite{N1}.

\medskip

In this paper, we address the homotopy theory of crossed module morphisms $\mathcal{X \to X'}$ in MCI which leads us to define an equivalence relation, therefore to construct a groupoid structure, without any restriction on $\mathcal{X}$ or $\mathcal{X'}$. This case should represents an undiscovered model category structure for the category of crossed modules, where all objects are both fibrant and cofibrant.

\medskip

As indicated in \cite{CL1}, we have the functorial relation between the categories of associative algebras, Leibniz algebras, Lie algebras, dialgebras and that of crossed modules in these categories which can be pictured as:
\begin{align}\label{functors}
\xymatrix@R=50pt@C=50pt{
	&               &       &       &\mathbf{DiAs} 
	\ar@{}[ddr]^(.1){}="a"^(.9){}="b" \ar_{Lb} "a";"b"
	\ar@{}[ddl]^(.1){}="a"^(.9){}="b" \ar^{As} "a";"b"   &       \\
	&\mathbf{XDiAs} 
	\ar@{}[ddr]^(.1){}="a"^(.9){}="b" \ar_{XLb} "a";"b"
	\ar@{}[ddl]^(.1){}="a"^(.9){}="b" \ar^{XAs} "a";"b"
	\ar@{<--}[urrr]^{J_i} &        &       &                   &       \\
	&               &       &\mathbf{As} 
	\ar@<0ex>[r]_{Lie_1}
	\ar@<2ex>@{}[ruu]^(.1){}="a"^(.9){}="b" \ar^{\subset} "a";"b"
	& \mathbf{Lie} \ar@<-1.5ex>[l]_{U}  \ar@<0ex>[r]_{\subset} & 
	\mathbf{Lb} \ar@<-1.5ex>[l]_{Lie_2}  \ar@<-2ex>@{}[luu]^(.1){}="a"^(.9){}="b" \ar_{U_d} "a";"b" \\
	\mathbf{XAs} \ar@<-1.5ex>[r]_{XLie_1}
	\ar@{<--}[urrr]^{I_i}  
	\ar@<2ex>@{}[ruu]^(.1){}="a"^(.9){}="b" \ar^{\subset} "a";"b"
	& \mathbf{XLie} \ar@<0ex>[l]_{XU} \ar@<-1.5ex>[r]_{\subset} & 
	\mathbf{XLb} \ar@<0ex>[l]_{XLie_2} 
	\ar@<-2ex>@{}[luu]^(.1){}="a"^(.9){}="b" \ar_{U_d} "a";"b"
	 \ar@{<--}[urrr]_{J'_i} &   &   & \\}
\end{align}
where all faces are commutative. Under this aspect, we will handle these crossed module structures and define the homotopy of morphisms by considering the particular cases of homotopy definition of crossed modules of MCI. On the other hand, one can see that the homotopy definitions given in \cite{JFM1, IJK} can also obtained from our generalized homotopy definition.  Moreover we also see that the adjoint crossed module functors given in \eqref{functors} are preserving the homotopy relation.

\medskip

The main result of this paper is that, the functor $X_1$ in \eqref{equality} preserves the homotopy, furthermore the homotopy equivalence between simplicial objects and crossed modules in MCI. However this property can not be extended to a groupoid functor since the groupoid structure of simplicial homotopies has not been discovered yet, even for groups or algebras.

\section{Preliminaries}

We will recall the main definitions and the statements from \cite{BDCU} which will be used in sequel.

\subsection{Modified Categories of Interest}

Let $\mathbb{C}$ be a category of groups with a set of operations $\Omega$ and
with a set of identities $\mathbb{E}$, such that $\mathbb{E}$ includes the
group identities and the following conditions hold. If $\Omega_{i}$ is the set
of $i$-ary operations in $\Omega$, then:

\begin{enumerate}
	\item[(a)] $\Omega=\Omega_{0}\cup\Omega_{1}\cup\Omega_{2}$;
	
	\item[(b)] the group operations (written additively : $0,-,+$) are elements of
	$\Omega_{0}$, $\Omega_{1}$ and $\Omega_{2}$ respectively. Let $\Omega
	_{2}^{\prime}=\Omega_{2}\setminus\{+\}$, $\Omega_{1}^{\prime}=\Omega
	_{1}\setminus\{-\}.$ Assume that if $\ast\in\Omega_{2}$, then $\Omega
	_{2}^{\prime}$ contains $\ast^{\circ}$ defined by $x\ast^{\circ}y=y\ast x$ and
	assume $\Omega_{0}=\{0\}$;
	
	\item[(c)] for each $\ast\in\Omega_{2}^{\prime}$, $\mathbb{E}$ includes the
	identity $x\ast(y+z)=x\ast y+x\ast z$;
	
	\item[(d)] for each $\omega\in\Omega_{1}^{\prime}$ and $\ast\in\Omega
	_{2}^{\prime}$, $\mathbb{E}$ includes the identities $\omega(x+y)=\omega
	(x)+\omega(y)$ and $\omega(x\ast y)=\omega(x)\ast\omega(y)$.
\end{enumerate}

\medskip

Let $C$ be an object of $\mathbb{C}$ and $x_{1},x_{2},x_{3}\in C$:

\bigskip

\noindent\textbf{Axiom 1.} For all $\ast\in\Omega_{2}^{\prime}$, we have:
\begin{align}\label{axiom}
x_{1}+(x_{2}\ast x_{3})=(x_{2}\ast x_{3})+x_{1}
\end{align}

\bigskip

\noindent\textbf{Axiom 2.} For each ordered pair $(\ast,\overline{\ast}%
)\in\Omega_{2}^{\prime}\times\Omega_{2}^{\prime}$ there is a word $W$ such
that:
\begin{gather*}
(x_{1}\ast x_{2})\overline{\ast}x_{3}=W(x_{1}(x_{2}x_{3}),x_{1}(x_{3}
x_{2}),(x_{2}x_{3})x_{1},\\
(x_{3}x_{2})x_{1},x_{2}(x_{1}x_{3}),x_{2}(x_{3}x_{1}),(x_{1}x_{3})x_{2}
,(x_{3}x_{1})x_{2}),
\end{gather*}
where each juxtaposition represents an operation in $\Omega_{2}^{\prime}$.

\bigskip

\begin{defn}
	\label{catint2}A category of groups with operations $\mathbb{C}$ satisfying
	conditions $(a)-(d)$, Axiom 1 and Axiom 2, will be called a modified
		category of interest (MCI). 
\end{defn}

As indicated in \cite{BDCU} the difference of this definition from the
original one of category of interest is the identity $\omega(x)\ast
\omega(y)=\omega(x\ast y),$ which is $\omega(x)\ast y=\omega(x\ast y)$ in the
definition of category of interest.

\begin{ex}
The categories of (pre)cat$^1$ objects in the category of Leibniz (Lie, Associative) algebras and dialgebras are all MCI which are not categories of interest.
\end{ex}

\begin{defn}
	Let $A$, $B\in\mathbb{C}$. An extension of $B$ by $A$ is a sequence:
	\begin{align*}
	\xymatrix{0\ar[r]&A\ar[r]^-{i}&E\ar[r]^-{p}&B\ar[r]&0}
	\end{align*}
	in which $p$ is surjective and $i$ is the kernel of $p$. We say that an
	extension is split if there is a morphism $s:B\to E$ such that
	$p\circ s=1_{B}$.
\end{defn}

\begin{defn}
	Suppose that $A,B$ is the objects of $\C$. We say that we have a set of
	actions of $B$ on $A$ if there is a map: $$f_{\ast}:A\times
	B\to A$$
	for all $\ast\in\Omega_{2}.$ A split extension
	of $B$ by $A,$ induces an action of $B$ on $A$ corresponding to the operations
	in $\mathbb{C}$ as the following:
	\begin{align*}
	b\cdot a & =s(b)+a-s(b), \\
	b\ast a & =s(b)\ast a,
	\end{align*}
	for all $b\in B$, $a\in A$ and $\ast\in\Omega_{2}{}^{\prime}.$ These actions will be called derived actions of $B$ on $A$. Alternatively, we can use also use the notation $\l$ to denote the actions which will be used in application section to make the different actions clear.
\end{defn}

\begin{defn}	\label{semidirect}
    Given an action of $B$ on $A,$ a semi-direct product
	$A\rtimes B$ is a universal algebra, whose underlying set is $A\times B$ and
	the operations are defined by
	\begin{align}\label{semiprop}
	\begin{split}
	\omega(a,b) & = (\omega\left(  a\right)  ,\omega\left(  b\right)  ),\\
	(a^{\prime},b^{\prime})+(a,b) & = (a^{\prime}+b^{\prime}\cdot a,b^{\prime
	}+b),\\
	(a^{\prime},b^{\prime})\ast(a,b) & = (a^{\prime}\ast a+a^{\prime}\ast
	b+b^{\prime}\ast a,b^{\prime}\ast b),
	\end{split}
	\end{align}
	for all $a,a^{\prime}\in A$ and $b,b^{\prime}\in B$.
\end{defn}

\begin{thm}
An action of $B$ on $A$ is a derived action if and only if
	$A\rtimes B$ is an object of $\mathbb{C}$.
\end{thm}

\begin{prop}
	 A set of actions of $B$ on $A$ in $\mathbb{C}_{G}$ is a set of
	derived actions \cite{BDCU} if and only if it satisfies the following conditions: 
	
	\begin{enumerate}

		\item[\textit{1.}] $0\cdot a=a$, 
		
		\item[\textit{2.}] $b\cdot(a_{1}+a_{2})=b\cdot a_{1}+b\cdot a_{2}$, 
		
		\item[\textit{3.}] $(b_{1}+b_{2})\cdot a=b_{1}\cdot(b_{2}\cdot a)$, 
		
		\item[\textit{4.}] $b\ast(a_{1}+a_{2})=b\ast a_{1}+b\ast a_{2}$, 
		
		\item[\textit{5.}] $(b_{1}+b_{2})\ast a=b_{1}\ast a+b_{2}\ast a$, 
		
		\item[\textit{6.}] $(b_{1}\ast b_{2})\cdot(a_{1}\ast a_{2})=a_{1}\ast a_{2}$, 
		
		\item[\textit{7.}] $(b_{1}\ast b_{2})\cdot(a\ast b)=a\ast b$, 
		
		\item[\textit{8.}] $a_{1}\ast(b\cdot a_{2})=a_{1}\ast a_{2}$, 
		
		\item[\textit{9.}] $b\ast(b_{1}\cdot a)=b\ast a$, 
		
		\item[\textit{10.}] $\omega(b\cdot a)=\omega(b)\cdot\omega(a)$, 
		
		\item[\textit{11.}] $\omega(a\ast b)=\omega(a)\ast\omega(b)$, \ 
		
		\item[\textit{12.}] $x\ast y+z\ast t=z\ast t+x\ast y$, 
	\end{enumerate}
	for all $\omega\in\Omega_{1}^{\prime}$, $\ast\in\Omega_{2}{}^{\prime}$, $b$,
	$b_{1}$, $b_{2}\in B$, $a,a_{1},a_{2}\in A$ and  $x,y,z,t\in A\cup B$
	whenever each side of $12$ has a sense.
\end{prop}

\subsection{Crossed Modules in MCI}

In the rest of the paper, $\C$ will denote an arbitrary MCI.

\begin{defn} \label{xmod}
	A crossed module in $\C$ given by a morphism $\d \colon E \to R$ together with a derived action of $R$ on $E$, such that the following relations called \textquotedblleft {Peiffer-Whitehead} relations\textquotedblright, hold:	
		\begin{description}
			\item[XM1)] $\d (r\cdot e)=r+\d(e)-r$ \quad and \quad $\d (r \ast e)=r \ast \d(e)$
			\item[XM2)] $\d (e)\cdot e'= e +e'-e$\quad and \quad $\d (e) \ast e'=e \ast e'$
		\end{description}
	for all $e,e' \in E$, $r \in R$ and $\ast\in\Omega_{2}'$. 
	
	\medskip
	
	\noindent Without the second relation we call it a precrossed module.
	
	\medskip
	
	From now on, $\X$ will denote a crossed module in $\C$ with being $\X=(E,R,\d)$.
\end{defn}

\begin{defn}
	Let $\X,\X'$ be two crossed modules. A crossed module morphism $f \colon \X \to \X'$ is a pair $f=(f_{1}\colon E \to E',f_{0}\colon R \to R')$ of morphisms in $\C$, making the diagram:
	\begin{equation}\label{commutes}
	\xymatrix@R=30pt@C=30pt{
		E
		\ar[r]^{\partial}
		\ar[d]_{f_{1}}
		& R
		\ar[d]^{f_{0}}
		\\
		E^{\prime}
		\ar[r]_{{\partial}^\prime}
		& R^{\prime}}\,
	\end{equation}
	commutative and also preserving the derived action of $R$ on $E$, means (for all $e \in E$ and $r \in R$):
	\begin{align*}
		f_{1}(r \cdot e)&=f_{0}(r)\cdot f_{1}(e) \\
		f_{1}(r \ast e)&=f_{0}(r)\ast f_{1}(e).
	\end{align*}
\end{defn}

	Consequently we have a category \textbf{XMod($\C$)} in MCI.

\subsection{Simplicial Objects in MCI}

We recall some simplicial data from \cite{PM2, GJ1, KP1}. 

\begin{defn}
	A simplicial object in $\C$ is a functor: $\Delta^{op} \to \mathcal{C}$ where $\Delta$ is the simplicial indexing category.
\end{defn}

An alternative definition of the simplicial object as the following:

\begin{defn}
	A simplicial object $\A$ in the category $\C$ is a collection of $\left\{ A_{n} \colon A_n \in Ob(\mathbb{C}) \, , \, n \in \mathbb{N} \right\}$ together with morphisms:
	\begin{equation*}
	\begin{tabular}{lllll}
	$d_{i}^{n-1}$ & $:$ & $A_{n}\longrightarrow A_{n-1}$ & $,$ & $0\leq i\leq
	n-1 $ \\
	$s_{j}^{n}$ & $:$ & $A_{n}\longrightarrow A_{n+1}$ & $,$ & $0\leq j\leq n$%
	\end{tabular}%
	\end{equation*}
	which are called face and degeneracies respectively (we will not use the superscripts in the calculations). 
	
	\medskip
	
	These homomorphisms are required to satisfy the following axioms,
	called simplicial identities:
	\begin{equation}\label{simpident}
	\begin{tabular}{llll}
	(i) & $d_{i}d_{j}=d_{j-1}d_{i}$ & if & $i<j$ \\
	(ii) & $s_{i}s_{j}=s_{j+1}s_{i}$ & if & $i\leq j$ \\
	(iii) & $d_{i}s_{j}=s_{j-1}d_{i}$ & if & $i<j$ \\
	& $d_{j}s_{j}=d_{j+1}s_{j}=id$ &  &  \\
	& $d_{i}s_{j}=s_{j}d_{i-1}$ & if & $i>j+1$ 
	\end{tabular}%
	\end{equation}
	
	Any simplicial object could be pictured as:
	\begin{align}\label{simplicial}
	\xymatrix@C=50pt{
		\mathbf{A} \: \doteq \: \: \ar@{.}[r] &
		{A}_{3}
		\ar@<2.25ex>[r] \ar@{.>}@<1.5ex>[r] \ar@{.>}@<0.75ex>[r] \ar@<0ex>[r] &
		{A}_{2}\ar@<3ex>[r]|{d_2} \ar@<1.5ex>[r]|{d_1} \ar@<0ex>[r]|{d_0}
		\ar@/^1pc/[l] \ar@{.>}@/^1.5pc/[l] \ar@/^2pc/[l]&
		{A}_{1}\ar@<1.5ex>[r]|{d_1} \ar[r]|{d_0}
		\ar@/^1pc/[l]|{s_0} \ar@/^1.5pc/[l]|{s_1} &
		{A}_{0}\ar@/^1pc/[l]|{s_0} }
	\end{align}
\end{defn}

\bigskip

\begin{defn}
	A simplicial map $f \colon \A \to \B$
	is a set of morphisms $f_{n}:A_{n} \to
	B_{n}$ commuting with all the face and degeneracy operators such that:
	\begin{align*}
	f_{q}d_{i} &=d_{i}f_{q+1} \\
	f_{q}s_{i} &=s_{i}f_{q-1}
	\end{align*}
	with the diagram:
	$$
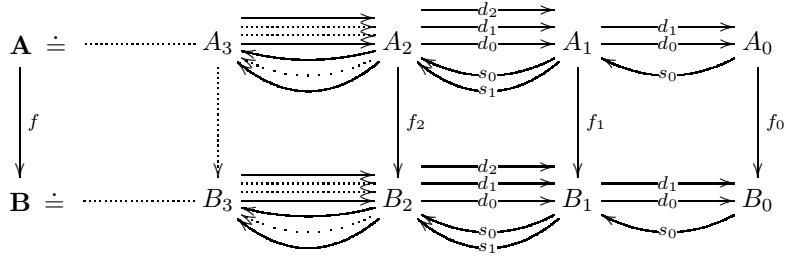
\xymatrix@R=45pt@C=50pt{
		*+[r]{\A \: \doteq \: \:} \ar@{.}[r]
		\ar@{}[d]^(.10){}="a"^(.90){}="b" \ar@[red] "a";"b"^{f}&
		{A}_{3}
		\ar@<2.25ex>[r] \ar@{.>}@<1.5ex>[r] \ar@{.>}@<0.75ex>[r] \ar@<0ex>[r]
		\ar@{}[d]^(.10){}="a"^(.90){}="b" \ar@[red]@{.>} "a";"b"&
		{A}_{2}\ar@<3ex>[r]|{d_2} \ar@<1.5ex>[r]|{d_1} \ar@<0ex>[r]|{d_0}
		\ar@/^0.5pc/[l] \ar@{.>}@/^1pc/[l] \ar@/^1.5pc/[l]
		\ar@{}[d]^(.10){}="a"^(.90){}="b" \ar@[red] "a";"b"^{f_{2}}&
		{A}_{1}\ar@<1.5ex>[r]|{d_1} \ar[r]|{d_0}
		\ar@/^1pc/[l]|{s_0} \ar@/^1.5pc/[l]|{s_1}
		\ar@{}[d]^(.10){}="a"^(.90){}="b" \ar@[red] "a";"b"^{f_{1}}&
		{A}_{0}\ar@/^1pc/[l]|{s_0}
		\ar@{}[d]^(.10){}="a"^(.90){}="b" \ar@[red] "a";"b"^{f_{0}}
		\\
		*+[r]{\B \: \doteq \: \:} \ar@{.}[r] &
		{B}_{3}
		\ar@<2.25ex>[r] \ar@{.>}@<1.5ex>[r] \ar@{.>}@<0.75ex>[r] \ar@<0ex>[r] &
		{B}_{2}\ar@<3ex>[r]|{d_2} \ar@<1.5ex>[r]|{d_1} \ar@<0ex>[r]|{d_0}
		\ar@/^0.5pc/[l] \ar@{.>}@/^1pc/[l] \ar@/^1.5pc/[l]&
		{B}_{1}\ar@<1.5ex>[r]|{d_1} \ar[r]|{d_0}
		\ar@/^1pc/[l]|{s_0} \ar@/^1.5pc/[l]|{s_1} &
		{B}_{0}\ar@/^1pc/[l]|{s_0}
	}$$ 
\end{defn}

	Consequently, we have thus defined the category of simplicial objects, which will denoted by \textbf{Simp($\mathcal{C}$)}.

\begin{defn}
	An $n$-truncated simplicial object is a simplicial object with objects  $A_{i}$  $(i \leq n)$. Therefore we can get a full subcategory of \textbf{Simp($\mathcal{C}$)}.
\end{defn}

\begin{defn}
	Given a simplicial object $\A$, the Moore Complex $(NA,\partial )$ of $\mathcal{A}$  is the chain complex defined by:
	\begin{equation*}
	NA_{n}=\underset{i=0}{\overset{n-1}{\mathlarger{\mathlarger{\mathlarger{{\cap}}}}}}Ker(d_{i}^{n})
	\end{equation*}%
	with the morphisms $\partial _{n} \colon NA_{n}\to NA_{n-1}$ induced from $%
	d_{n}^{n-1} $ by restriction.
\end{defn}

\begin{defn}
	Let $(NA,\partial )$ be a Moore complex of a simplicial object $\A$. We call this Moore Complex with length $n$, iff $NA_{i}$ is equal to $\{0\}$, for each	$i>n$. We denote the category of simplicial objects with Moore Complex of length $n$ by \textbf{Simp$ \, (\mathcal{C})_{\leq n}$}.
\end{defn}

\subsubsection{Simplicial Homotopy}

\begin{defn}
	Let $f,g \colon \A \to \B$ be simplicial maps. If there exist the family of morphisms of $\mathcal{C}$ defined as $h_{i}^{n}:A_{n}\longrightarrow B_{n+1}$, $0\leq
	i\leq q$ which satisfies:
	\begin{equation}\label{homident}
	\begin{tabular}{llll}
	(i) & $d_{0}h_{0}=f\ \ ,\ \ d_{q+1}h_{q}=g$ &  &  \\
	(ii) & $d_{i}h_{j}=h_{j-1}d_{i}$ & if & $i<j$ \\
	& $d_{j+1}h_{j+1}=d_{j+1}h_{j}$ &  &  \\
	& $d_{i}h_{j}=h_{j}d_{i-1}$ & if & $i>j+1$ \\
	(iii) & $s_{i}h_{j}=h_{j+1}s_{i}$ & if & $i\leq j$ \\
	& $s_{i}h_{j}=h_{j}s_{i-1}$ & if & $i>j$ \\
	\end{tabular}
	\end{equation}
	then we say that the collection of $\left\{ h_i \right\}$ defines a homotopy \cite{PM2} connecting $f$ to $g$ and denote it $f\simeq g$. All fits in the diagram:
	$$
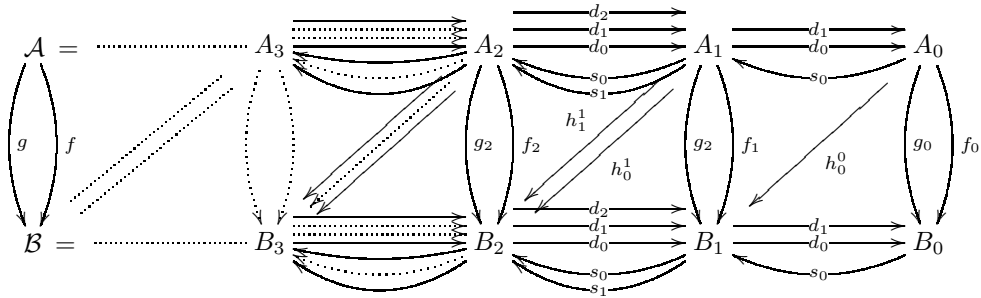
\xymatrix@R=60pt@C=65pt{
		*+[r]{\mathcal{A} \: = \: \:} \ar@{.}[r]
		\ar@/^0.75pc/[d]^{f}
		\ar@/_0.75pc/[d]^{g}&
		{A}_{3}
		\ar@<2.25ex>[r] \ar@{.>}@<1.5ex>[r] \ar@{.>}@<0.75ex>[r] \ar@<0ex>[r]
		\ar@{.>}@/^0.75pc/[d]
		\ar@{.>}@/_0.75pc/[d]
		\ar@{}[dl]^(.15){}="a"^(.85){}="b" \ar@[red]@{.}@<0.7ex> "a";"b"
		\ar@{}[dl]^(.15){}="a"^(.85){}="b" \ar@[red]@{.}@<-0.7ex> "a";"b"&
		{A}_{2}\ar@<3ex>[r]|{d_2} \ar@<1.5ex>[r]|{d_1} \ar@<0ex>[r]|{d_0}
		\ar@/^0.5pc/[l] \ar@{.>}@/^1pc/[l] \ar@/^1.5pc/[l]
		\ar@/^0.75pc/[d]^{f_{2}}
		\ar@/_0.75pc/[d]^{g_{2}}
		\ar@{}[dl]^(.15){}="a"^(.85){}="b" \ar@[red]@<0.9ex> "a";"b"
		\ar@{}[dl]^(.15){}="a"^(.85){}="b" \ar@[red]@{.>} "a";"b"
		\ar@{}[dl]^(.15){}="a"^(.85){}="b" \ar@[red]@<-0.9ex> "a";"b"&
		{A}_{1}\ar@<1.5ex>[r]|{d_1} \ar[r]|{d_0}
		\ar@/^1pc/[l]|{s_0} \ar@/^1.5pc/[l]|{s_1}
		\ar@/^0.75pc/[d]^{f_{1}}
		\ar@/_0.75pc/[d]^{g_{2}}
		\ar@{}[dl]^(.15){}="a"^(.85){}="b" \ar@[red]@<0.7ex> "a";"b"^{{h_{0}^{1}}}
		\ar@{}[dl]^(.15){}="a"^(.85){}="b" \ar@[red]@<-0.7ex> "a";"b"_{{h_{1}^{1}}}&
		{A}_{0}\ar@/^1pc/[l]|{s_0}
		\ar@{}[dl]^(.15){}="a"^(.85){}="b" \ar@[red] "a";"b"^{h_{0}^{0}}
		\ar@/^0.75pc/[d]^{f_{0}}
		\ar@/_0.75pc/[d]^{g_{0}}
		\\
		*+[r]{\mathcal{B} \: = \: \:} \ar@{.}[r] &
		{B}_{3}
		\ar@<2.25ex>[r] \ar@{.>}@<1.5ex>[r] \ar@{.>}@<0.75ex>[r] \ar@<0ex>[r] &
		{B}_{2}\ar@<3ex>[r]|{d_2} \ar@<1.5ex>[r]|{d_1} \ar@<0ex>[r]|{d_0}
		\ar@/^0.5pc/[l] \ar@{.>}@/^1pc/[l] \ar@/^1.5pc/[l]&
		{B}_{1}\ar@<1.5ex>[r]|{d_1} \ar[r]|{d_0}
		\ar@/^1pc/[l]|{s_0} \ar@/^1.5pc/[l]|{s_1} &
		{B}_{0}\ar@/^1pc/[l]|{s_0}
	}$$ 
\end{defn}

\section{Homotopy of Crossed Modules in MCI}

In the rest of the paper, we fix two arbitrary crossed modules $\mathcal{X}=(E,R,\d)$ and $\mathcal{X'}=(E',R',\d')$ in $\C$.

\subsection{Derivation and Homotopy}

\begin{defn}\label{sderivation}
	Let $f_{0}\colon R\to R^{\prime }$ be a morphism in $\mathcal{C}$. An $%
	f_{0}$-derivation $s\colon R\to E^{\prime }$ is a map satisfying:
	\begin{align}\label{s-der}
	 \begin{split}
		s(g + h) & = \Big( f_0(-h) \cdot s(g) \Big) + s(h)  \\
		s(g \ast h) & = f_{0}(g) \ast s(h)+ f_{0}(h) \ast^{\circ} s(g) + s(g) \ast s(h)
	 \end{split}
	\end{align}
	for all $g,h \in R$.
\end{defn}

\begin{lem}\label{s-prop}
	If $s$ is an $f_0$-derivation, then:
	\begin{itemize}
		\item $s(0)=0$
		\item $s(-g)=f_0(g) \cdot (-s(g))$
		\item $s(g+h-g)=f_0(g) \cdot \Big( f_0(-h) \cdot s(g) + s(h) \Big) + s(-g)$
	\end{itemize}
\end{lem}

\begin{pf}
	Easy calculations.
\end{pf}

\begin{lem}
	Let $f \colon \X \to \X'$ be a crossed module morphism. Any $f_0$ derivation $s \colon R \to E'$ can be seen as a (unique) morphism in $\C$ with being:
	\begin{align*}
		\phi \colon r \in R \mapsto \big( f_0(r) , s(r) \big) \in R' \ltimes E'.
	\end{align*}
\end{lem}

\begin{pf}
	Direct checking by using \eqref{s-der} and \eqref{semiprop}.
\end{pf}

\begin{thm}
	Let $f \colon \X \to \X'$ be a crossed module morphism. If $s$ is an $f_{0}$-derivation, and if we define $%
	g=(g_{1},g_{0})$ as  (where $e\in E$ and $r\in R$):
	\begin{align}\label{g-defn}
	&g_{0}(r)=f_{0}(r)+(\partial ^{\prime } \circ s)(r), 
	&&g_{1}(e)=f_{1}(e)+(s \circ \partial )(e),
	\end{align}
	then $g$ is also defines a crossed module morphism $\mathcal{X \to X'}$.
\end{thm}
	
	\begin{pf} To make the formula more compact, in the rest of the paper, we do not use $\circ$ to denote composition in the proofs. Since for all $r,r' \in R$:
		\begin{align*}
		g_{0}(r \ast r') & =  f_{0}(r \ast r')+\d' s(r \ast r') \\
		& =  f_{0}(r) \ast f_{0}(r')+\d' \Big( f_{0}(r) \ast s(r')+f_0(r') \ast^{\circ} s(r)+s(r) \ast s(r') \Big) \\
		& =  f_{0}(r) \ast f_{0}(r')+\d' ( f_{0}(r) \ast s(r'))+\d' (f_0(r') \ast^{\circ} s(r))+\d' (s(r) \ast s(r')) \\
		& =  f_{0}(r) \ast f_{0}(r')+f_{0}(r) \ast \d' s(r')+f_{0}(r') \ast^{\circ} \d's(r)+\d's(r) \ast \d's(r') \\
		& =  (f_{0}(r)+\d's(r)) \ast (f_{0}(r')+\d's(r')) \\
		& =  g_{0}(r) \ast g_{0}(r')
		\end{align*}
		and
		    \begin{align*}
			g_{0}(r + r') & =  f_{0}(r + r')+\d' s(r + r') \\
			& =  f_{0}(r) + f_{0}(r')+\d' \Big( \big( f_0(-r') \cdot s(r) \big) + s(r') \Big) \\
			& = f_{0}(r) + f_{0}(r') + \d' \big( f_0(-r') \cdot s(r) \big) + \d's(r')  \\
			& = f_{0}(r) + f_{0}(r') - f_0(r') + \d's(r) + f_0(r') + \d's(r')  \\
			& =  f_{0}(r) + \d's(r) + f_0(r') + \d's(r')  \\
			& = g_0 (r) + g_0 (r')
			\end{align*}
			$g_0$ is a morphism in $C$; similarly $g_1$. It is also easy to check that the diagram \ref{commutes}	commutes. Finally; $g_1$ preserves the derived actions of $R$ on $E$. Indeed:
		\begin{align*}
		g_{1}(r \ast e) & =  f_{1}(r\ast e)+s\d (r\ast e) \\
		& =  f_{0}(r)\ast f_{1}(e)+s(r\d (e)) \\
		& =  f_{0}(r)\ast f_{1}(e)+f_{0}(r)\ast s\d (e)+ f_{0}\d (e) \ast^{\circ} s(r) +s(r) \ast s(\d (e)) \\
		& =  f_{0}(r)\ast f_{1}(e)+f_{0}(r)\ast s\d (e)+\d' f_{1}(e) \ast^{\circ} s(r) +s(r) \ast s(\d (e)) \\
		& = f_{0}(r)\ast f_{1}(e)+f_{0}(r)\ast s\d(e)+s(r) \ast f_{1}(e)+s(r) \ast s(\d (e)) \\
		& =  f_{0}(r)\ast f_{1}(e)+f_{0}(r)\ast s\d
		(e)+\d's(r)\ast f_{1}(e)+\d's(r)\ast s\d (e) \\
		& =  (f_{0}(r)+\d's(r))\ast(f_{1}(e)+s\d (e)) \\
		& =  g_{0}(r) \ast g_{1}(e)
		\end{align*}
		for all $r \in R $ and $e \in E$. On the other hand:
		\begin{align*}
			g_{o}(r) \cdot g_1(e) & = \Big( f_0(r) + \d's(r) \Big) \cdot \Big( f_1(e) + s\d(e) \Big) \\
			& = \Big( f_0(r) + \d's(r) \Big) \cdot f_1(e) + \Big( f_0(r) + \d's(r) \Big) \cdot s\d(e) \\
			& = f_0(r) \cdot \Big( \d's(r) \cdot f_1(e) \Big) + f_0(r) \cdot \Big( \d's(r) \cdot s\d(e) \Big) \\
			& = f_0(r) \cdot \Big( s(r) +  f_1(e) -s(r) \Big) + f_0(r) \cdot \Big( s(r) + s\d(e) - s(r) \Big) \\
			& = f_0(r) \cdot s(r) + f_0(r) \cdot f_1(e) + f_0(r) \cdot (-s(r)) + f_0(r) \cdot s(r) + f_0(r)  \cdot s\d(e) + f_0(r) \cdot (-s(r)) \\
			& = f_0(r) \cdot s(r) + f_0(r) \cdot f_1(e) + f_0(r)  \cdot s\d(e) + f_0(r) \cdot (-s(r)) \\
		\end{align*}
		and also:
		\begin{align*}
			g_1(r \cdot e) & = f_1(r \cdot e) + s\d (r \cdot e) \\
			& = f_0 (r) \cdot f_1 (e) + s \big(r + \d(e) - r  \big) \\
			& = f_0 (r) \cdot f_1 (e) + f_0(r) \cdot \Big( f_0(\d(e)) \cdot s(r) + s(\d(e)) \Big) + s(-r) \\
			& = f_0 (r) \cdot f_1 (e) + f_0(r) \cdot \Big( \d'(f_1(-e)) \cdot s(r) + s(\d(e)) \Big) + s(-r) \\
			& = f_0 (r) \cdot f_1 (e) + f_0(r) \cdot \Big( f_1(-e) + s(r) + f_1(e) + s(\d(e)) \Big) +  \\
			& = f_0(r) \cdot s(r) + f_0(r) \cdot f_1(e) + f_0(r)  \cdot s\d(e) + s(-r).
		\end{align*}
		for all $r \in R$ and $e \in E$ and by using Lemma \ref{s-prop}:
		\begin{align*}
		g_{1}(r \cdot e) =  g_{o}(r) \cdot g_1(e).
		\end{align*}
		Therefore $g=(g_{1},g_{0})$ is a crossed module morphism between $\X \to \X'$.
	\end{pf}
	
	\begin{defn}\label{homotopy}
		In the condition of the previous theorem, we write
		$f\xra{(f_{0},s)}g$ or shortly $f \simeq g$,
		and say that $(f_{0},s)$ is a \textbf{homotopy (or derivation)}  connecting $f$ to $g$. 
		
		\medskip
		
		As a consequence of this homotopy definition, we can give the following: 
		
		\medskip
		
		Let $\X, \X'$ be crossed modules. If there exist crossed module morphisms $f \colon \X \to \X'$ and $g \colon \X' \to \X$ such that $f \circ g \simeq id_{\X'}$ and $g \circ f \simeq id_{\X}$; we say that the crossed modules $\X$ and $\X'$ are \textbf{homotopy equivalent}, which denoted by $\X \simeq \X'$.
	\end{defn}
	
	\begin{rem}
	In the calculations above, we used  the second Peiffer-Whitehead relation, from Definition \ref{xmod}. So that this homotopy definition does not hold for precrossed modules (see details in \cite{JFM1} for the group case).
	\end{rem}

\subsection{A Groupoid}

Now we construct a groupoid structure which is induced from homotopy of crossed module morphisms in $\C$.

\begin{lem}[Identity]
	Let $f=(f_1,f_0)$ be a crossed module morphism $\X \to \X'$. The null function $
	0_{s}: r \in R  \longmapsto 0_{E^{\prime }} \in   E^{\prime } 
	$
	defines an $f_{0}$-derivation connecting $f$ to $f$.
\end{lem}
	
	\begin{pf}
    Easy calculations.
	\end{pf}

\begin{lem}[Inverse]
	Let $f=(f_1,f_0)$ and $g=(g_1,g_0)$ be crossed module morphisms $\X \to X'$ and $s$ be an $f_{0}$-derivation connecting $%
	f $ to $g$. Then, the map $\bar{s}=-s\colon R \to E'$, with 
	$\bar{s}(r)=-s(r)$, where $r \in R$,
	is a $g_{0}$-derivation connecting $g$ to $f$.
\end{lem}
	
	\begin{pf}
		Since $s$ is an $f_{0}$-derivation connecting $f$ to $g$, we have (for all $r,r' \in R$):
		\begin{align*}
		&f_{0}(r)=g_{0}(r)+(\partial ^{\prime }\circ \bar{s})(r), &\textrm{ and } &&
		f_{1}(r)=g_{1}(r)+(\bar{s}\circ \partial )(r).
		\end{align*}
		Moreover $\bar{s}$ is a $g_{0}$-derivation, since:
		\begin{align*}
		\bar{s}(g+h) &= - \big(s(g + h) \big) \\
		&= -\Big( \big( f_0(-h) \cdot s(g) \big) + s(h) \Big) \\
		&= - s(h) - \Big( f_0(-h) \cdot s(g) \Big) \\
		&=  - s(h) + \Big( f_0(-h) \cdot (-s(g)) \Big) \\
		& = f_0(-h) \cdot s(-h) + f_0(-h) \cdot \bar{s}(g) + f_0(-h) \cdot (-s(-h))  + \bar{s}(h) \\
		&= f_0(-h) \cdot \Big(s(-h) + \bar{s}(g)-s(-h) \Big) + \bar{s}(h)  \\
		&= f_0 (-h) \cdot \Big( \d's(-h) \cdot \bar{s}(g) \Big) + \bar{s}(h) \\
		&= \Big( \big( f_0(-h) + \d's(-h) \big) \cdot \bar{s}(g) \Big) + \bar{s}(h) \\
		&= \Big( g_0(-h) \cdot \bar{s}(g) \Big) + \bar{s}(h).
		\end{align*} 
		also:
				\begin{align*}
				\bar{s}(g \ast h) &= - \big(s(g \ast h) \big) \\
				&= -\Big(f_{0}(g) \ast s(h)+ f_{0}(h) \ast^{\circ} s(g) + s(g) \ast s(h) \Big) \\
				&= - f_{0}(g) \ast s(h) - f_{0}(h) \ast^{\circ} s(g) - s(g) \ast s(h) + s(g) \ast s(h) - s(g) \ast s(h) \\
				&= - f_{0}(g) \ast s(h) - s(g) \ast s(h)  - f_{0}(h) \ast^{\circ} s(g) - s(h) \ast^{\circ} s(g) + s(g) \ast s(h) \\
				& = - f_{0}(g) \ast s(h) - \d's(g) \ast s(h)  - f_{0}(h) \ast^{\circ} s(g) - \d's(h) \ast^{\circ} s(g) + s(g) \ast s(h) \\
				& = f_{0}(g) \ast (-s(h)) + \d's(g) \ast (-s(h))  + f_{0}(h) \ast^{\circ} (-s(g)) - \d's(h) \ast^{\circ} (-s(g)) + s(g) \ast s(h) \\
				& = \big( f_0(g) + \d's(g) \big) \ast (-s(h)) + \big( f_0(h) + \d's(h) \big) \ast (-s(g)) + s(g) \ast s(h) \\
				&= g_{0}(g) \ast (-s(h)) + g_{0}(h) \ast^{\circ} (-s(g))  + s(g) \ast s(h) \\
				&= g_{0}(g) \ast \bar{s}(h) + g_{0}(h) \ast^{\circ} \bar{s}(g) + \bar{s}(g) \ast \bar{s}(h).
				\end{align*} 
				for all $g,h \in R$. Note that in the first part of the proof we frequently used Lemma \ref{s-prop}
	\end{pf}

\begin{lem}[Concatenation]
	Let $f=(f_1 , f_0), \ g=(g_1 , g_0)$ and $k=(k_1 , k_0)$ be crossed module
	morphisms $\X \to \X'$, $s$ be an $f_{0}$%
	-derivation connecting $f$ to $g$, and $s^{\prime }$ be a $g_{0}$%
	-derivation connecting $g$ to $k$. Then the linear map $(s+s')\colon R \to E'$, such that $(s+s^{\prime })(r)=s(r)+s^{\prime }(r)$, defines an $f_{0}$-derivation (therefore a homotopy) connecting $f$ to $k$.
	
	\medskip
	
	\begin{pf}
		We know that $f\xra{ (f_0,s) }g$ and $g\xra{(g_0,s')} k$. Therefore, by definition:
		\begin{align*}
		&k_{0}(r)=f_{0}(r)+(\partial ^{\prime }\circ (s+s^{\prime }))(r), && k_{1}(e)=f_{1}(e)+((s+s^{\prime })\circ \partial )(e).
		\end{align*}%
		Since:
			\begin{align*}
			(s+s')(g+h) & = s(g+h)+s'(g+h) \\
            & = \Big( f_0(-h) \cdot s(g) \Big) + s(h) + \Big( g_0(-h) \cdot s'(g) \Big) + s'(h) \\
            & = \Big( f_0(-h) \cdot s(g) \Big) + s(h) + \Big( \big( f_0(-h) + \d's(-h) \big) \cdot s'(g) \Big) + s'(h) \\
            & = \Big( f_0(-h) \cdot s(g) \Big) + s(h) + \Big( f_0(-h) \cdot \big( \d's(-h) \cdot s'(g) \big) \Big) + s'(h) \\
            & = \Big( f_0(-h) \cdot s(g) \Big) + s(h) + \Big( f_0(-h) \cdot \big( s(-h) + s'(g) - s(-h) \big) \Big) + s'(h) \\
            & = \Big( f_0(-h) \cdot s(g) \Big) + s(h) + f_0(-h) \cdot \Big( f_0(h) \cdot (-s(h)) \Big)+ f_0(-h) \cdot s'(g) - f_0(-h) \cdot \Big( f_0(h) \cdot (-s(h)) \Big) + s'(h) \\
            & = f_0(-h) \cdot s(g) + f_0(-h) \cdot s'(g) + s(h) + s'(h) \\
            & = \Big( f_0 (-h) \cdot (s+s')(g) \Big) (s+s')(h)
			\end{align*}
			and also:
		\begin{align*}
		(s+s')(g \ast h) & = s(g \ast h)+s'(g \ast h)
		\\
		& = f_{0}(g)  \ast  s(h)+ f_{0}(h)  \ast^{\circ} s(g) +s(g)  \ast  s(h)+g_{0}(g)  \ast 
		s'(h) + g_{0}(h)  \ast^{\circ} s'(g) +s'(g) \ast s'(h) \\
		& = f_{0}(g)  \ast  s(h)+ f_{0}(h)  \ast^{\circ} s(g) +s(g) \ast s(h) + \big(f_0(g) + (\d' \circ s) (g) \big) \ast s'(h) \\
		& \quad + \big(f_0(h) + (\d' \circ s) (h) \big) \ast^{\circ} s'(g) +s'(g) \ast s'(h)  \\
		& = f_{0}(g)  \ast  s(h)+ f_{0}(h)  \ast^{\circ} s(g) +s(g) \ast s(h) + f_0(g) \ast s'(h) +  (\d' \circ s) (g) \ast s'(h) \\
		& \quad + f_0(h) \ast^{\circ} s'(g) + (\d' \circ s) (h) \ast^{\circ} s'(g) +s'(g) \ast s'(h)  \\
		& = f_{0}(g)  \ast  s(h)+ f_{0}(h)  \ast^{\circ} s(g) +s(g) \ast s(h) + f_0(g) \ast s'(h) +  s (g) \ast s'(h) \\
		& \quad + f_0(h) \ast^{\circ} s'(g) + s (h) \ast^{\circ} s'(g) +s'(g) \ast s'(h) \\ 
		& = f_{0}(g)  \ast  (s+s')(h)+ f_0(h)  \ast^{\circ} (s+s')(g) +(s+s')(g)  \ast (s+s')(h)
		\end{align*}
		for all $g,h \in R$; $(s+s')$ is an $f_{0}$-derivation connecting $f$ to $k$.
	\end{pf}
\end{lem}

\begin{rem}
	Notice that, in the proofs of previous two lemmas, we frequently used the property \eqref{axiom} and the crossed module axioms given in Definition \ref*{xmod}.
\end{rem}

Now we can give the following:

\begin{cor}{Let $\mathcal{X}, \mathcal{X^{\prime}}$ be two arbitrary but fixed crossed modules in $\C$. We have a groupoid ${\rm HOM}(\mathcal{X}, \mathcal{X'})$,  whose objects are the crossed module morphisms $\mathcal{%
X \rightarrow X'}$,  the  morphisms being their homotopies. In particular the  relation below, for crossed module morphisms $\mathcal{X} \to \mathcal{X'}$, is an equivalence relation:}
	\begin{equation*}
	\text{\textquotedblleft }f\simeq g\Longleftrightarrow \text{there exists an }%
	f_{0}\text{-derivation }s \text{ connecting } f \text{ with } g\text{\textquotedblright }.
	\end{equation*}
\end{cor}

\begin{pf}
	Follows from previous three lemmas.
\end{pf}

\section{From Simplicial Homotopy to Crossed Module Homotopy}

It is a well-know equivalence that, for a (modified) category of interest $\C$; category of crossed modules are equivalent to category of simplicial objects with Moore complex of length one \cite{YO1} with the  functors:
    $$\xymatrix@C=15pt@R=30pt{
    	Simp(\C)_{\leq 1} \ar[rr]^{X_1} & & XMod \ar[dl] \\
    	& Tr_{1}Simp(\C)_{\leq 1} \ar[ul]^{t_1}}
    $$ 
\noindent In this section, we will enrich the functor $X_1 \colon Simp(\C)_{\leq 1} \to XMod$ by exploring its relation with homotopy.

\medskip

Now let us recall how the functor $X_1$ works:

\medskip

Suppose that $\A$ is a simplicial algebra with Moore complex of length one, as seen on \eqref{simplicial}. We can construct a crossed module by the functor $X_1$ as the following:
	\begin{enumerate}
		\item Put $R=NA_{0}=A_{0}$ and$\ E=NA_{1}=Ker(d_0)$
		
		\item $R$ act on $E$ by:
		   \begin{align*}
		   	r \cdot e & = s_0 (r) + e - s_0 (r') \\
		   	r \ast e & = s_0 (r) \ast e
		   \end{align*}
		
		\item $\partial =d_{1}^{0}$ (restricted to $E$)
	\end{enumerate}
	
	Then we get the crossed module $(E,R,\d)$ with being:	
	\begin{align*}
	Ker(d_{0}^{0})\overset{d_{1}^{0}}{\longrightarrow }A_{0}
	\end{align*}

\begin{thm}\label{main1} 
	The functor $X_1$ preserves the homotopy. On other words, let $\A$ and $\B$ be any simplicial objects with Moore complex of length one and$\ f,g \colon \A \to \B$ are simplicial
		maps such that $h:f\simeq g$. Then: $$X_1(f) \simeq X_1(g).$$
	
\end{thm}

\begin{pf} 
	Let us define a map:
		\begin{align*}\label{transform}
		\zeta \colon \{h_i\} \longmapsto \left\{ -s_0f_0+h_0^0  \right\}
		\end{align*}
		where $h=\{h_i\}$ is the homotopy of simplicial maps $f\simeq g$. 
		
		\medskip
		
		Our claim is that: $\zeta$ defines a homotopy between crossed module morphisms: $$ X_1 (f) \xra{(f_0 , \zeta \{h_i\} )} X_1(g).$$
		
		Diagrammatically:
		
		$$
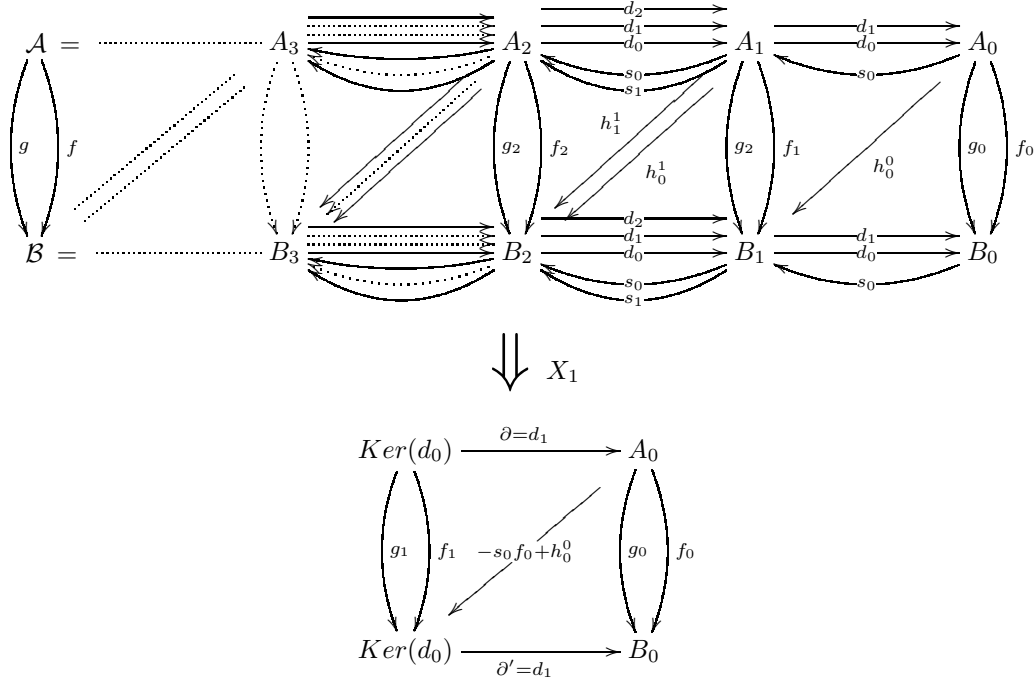
\xymatrix@R=65pt@C=70pt{
			*+[r]{\mathcal{A} \: = \: \:} \ar@{.}[r]
			\ar@/^0.75pc/[d]^{f}
			\ar@/_0.75pc/[d]^{g}&
			{A}_{3}
			\ar@<2.25ex>[r] \ar@{.>}@<1.5ex>[r] \ar@{.>}@<0.75ex>[r] \ar@<0ex>[r]
			\ar@{.>}@/^0.75pc/[d]
			\ar@{.>}@/_0.75pc/[d]
			\ar@{}[dl]^(.15){}="a"^(.85){}="b" \ar@[red]@{.}@<0.7ex> "a";"b"
			\ar@{}[dl]^(.15){}="a"^(.85){}="b" \ar@[red]@{.}@<-0.7ex> "a";"b"&
			{A}_{2}\ar@<3ex>[r]|{d_2} \ar@<1.5ex>[r]|{d_1} \ar@<0ex>[r]|{d_0}
			\ar@/^0.5pc/[l] \ar@{.>}@/^1pc/[l] \ar@/^1.5pc/[l]
			\ar@/^0.75pc/[d]^{f_{2}}
			\ar@/_0.75pc/[d]^{g_{2}}
			\ar@{}[dl]^(.15){}="a"^(.85){}="b" \ar@[red]@<0.9ex> "a";"b"
			\ar@{}[dl]^(.15){}="a"^(.85){}="b" \ar@[red]@{.>} "a";"b"
			\ar@{}[dl]^(.15){}="a"^(.85){}="b" \ar@[red]@<-0.9ex> "a";"b"&
			{A}_{1}\ar@<1.5ex>[r]|{d_1} \ar[r]|{d_0}
			\ar@/^1pc/[l]|{s_0} \ar@/^1.5pc/[l]|{s_1}
			\ar@/^0.75pc/[d]^{f_{1}}
			\ar@/_0.75pc/[d]^{g_{2}}
			\ar@{}[dl]^(.15){}="a"^(.85){}="b" \ar@[red]@<0.7ex> "a";"b"^{{h_{0}^{1}}}
			\ar@{}[dl]^(.15){}="a"^(.85){}="b" \ar@[red]@<-0.7ex> "a";"b"_{{h_{1}^{1}}}&
			{A}_{0}\ar@/^1pc/[l]|{s_0}
			\ar@{}[dl]^(.15){}="a"^(.85){}="b" \ar@[red] "a";"b"^{h_{0}^{0}}
			\ar@/^0.75pc/[d]^{f_{0}}
			\ar@/_0.75pc/[d]^{g_{0}}
			\\
			*+[r]{\mathcal{B} \: = \: \:} \ar@{.}[r] &
			{B}_{3}
			\ar@<2.25ex>[r] \ar@{.>}@<1.5ex>[r] \ar@{.>}@<0.75ex>[r] \ar@<0ex>[r] &
			{B}_{2}\ar@<3ex>[r]|{d_2} \ar@<1.5ex>[r]|{d_1} \ar@<0ex>[r]|{d_0}
			\ar@/^0.5pc/[l] \ar@{.>}@/^1pc/[l] \ar@/^1.5pc/[l]&
			{B}_{1}\ar@<1.5ex>[r]|{d_1} \ar[r]|{d_0}
			\ar@/^1pc/[l]|{s_0} \ar@/^1.5pc/[l]|{s_1} &
			{B}_{0}\ar@/^1pc/[l]|{s_0}
		}$$ 
		\begin{equation*}
		{\Huge
			\begin{tabular}[t]{r}
			$\Downarrow $%
			\end{tabular}%
		}{\Large X_1}
		\end{equation*}
		$$ \xymatrix@R=60pt@C=60pt{
			Ker(d_{0})
			\ar[r]^{\d=d_{1}}
			\ar@/^0.75pc/[d]^{f_{1}}
			\ar@/_0.75pc/[d]^{g_{1}}
			& A_{0}
			\ar@{}[dl]^(.15){}="a"^(.85){}="b" \ar@[red] "a";"b"|-{-s_{0}f_{0}+{h_{0}^{0}}}
			\ar@/^0.75pc/[d]^{f_{0}}
			\ar@/_0.75pc/[d]^{g_{0}}
			\\
			Ker(d_{0})
			\ar[r]_{\d'=d_{1}}
			& B_{0}
		}$$
		Recall the construction of $X_1(\A)$; therefore we have $X_1(f_0)=f_0$ and also $X_1(f_1)=f_1$ by its restriction \cite{YO1}.
		
		\medskip
	
	To reduce the calculations below, we will put $H(a)= \big(-s_{0}f_{0}+{h_{0}^{0}} \big) (a)$.
	
	\medskip
	
	\textbf{(i)} First of all the map is $\zeta$ well defined since (for all $a \in A_0$):
	\begin{align*}
	d _{0} \Big( \big(-s_{0}f_{0}+{h_{0}^{0}} \big) (a) \Big) & = d_{0} \Big(-(s_{0}f_{0})(a) + h_{0}^{0}(a) \Big) \\
	& =-d_{0} \big( (s_{0}f_{0})(a) \big) + d_{0}(h_{0}^{0}(a)) \\
	& =- (d_{0} s_{0})(f_{0}(a)) + f_{0}(a) \\
	& =-f_{0}(a)+f_{0}(a) \\
	& = 0_{B_{0}}
	\end{align*}
	which means:
	\begin{align*}
		Im(\zeta) \subseteq Ker(d_{0})=NB_{1}.
	\end{align*}
	
	\textbf{(ii)} Now we need to check the conditions given in \eqref{g-defn} for $H$. On other words the following conditions are to satisfy:
	\begin{align*}
	g_{0}(a) & = f_{0}(a) + (\d' \circ H)(a) \\
	g_{1}(e) & = f_{1}(e) + (H \circ \d)(e)
	\end{align*}
	
	It is clear that:
	\begin{align*}
	\d' \circ H & = d_{1} \Big( -s_{0}f_{0} + h_{0}^{0} \Big) \\
	& =  -d_{1}s_{0}f_{0} + d_{1}h_{0}^{0} \\
	& = -f_{0} + d_{1}h_{0}^{0} \\
	& = -f_{0} + g_{0}
	\end{align*}
	which leads to:
	\begin{align*}
		g_{0}(a) = f_{0}(a) + (\d' \circ H)(a)
	\end{align*}
	for all $a \in A_0$.
	
	\bigskip
	
	For the second condition required, we get:
	\begin{align}\label{continue}
	\begin{split}
	H \circ \d & = 	(-s _{0}f_{0} + h_{0}^{0})d_{1} \\
	&= -s	_{0}f_{0}d_{1} + h_{0}^{0}d_{1} \\
	& =  -s _{0}d _{0}h_{0}^{0}d_{1} + h_{0}^{0}d_{1} \\
	& = -s _{0}d _{0}d _{2}h_{0}^{1} + d _{2}h_{0}^{1} 	\\
	& =  -s _{0}d _{1}d _{0}h_{0}^{1} + d _{2}h_{0}^{1} \\
	& =  -d _{2}s _{0}d _{0}h_{0}^{1} + d _{2}h_{0}^{1}\\
	& =  -d _{2}s _{0}d	_{0}h_{0}^{1}+ d _{2}h_{0}^{1} -g_1 + f_1 -f_1 + g_1  \\
	& =  \big( -d _{2}s _{0}d_{0}h_{0}^{1} + d _{2}h_{0}^{1} -d _{2}h_{1}^{1}+d _{0}h_{0}^{1} \big)-f_1 + g_1 \\
	& =  d _{2} \Big(-s _{0}d_{0}h_{0}^{1}+h_{0}^{1}-h_{1}^{1}+s _{1}d _{0}h_{0}^{1} \Big)-f_1 + g_1 \\
	\end{split}
	\end{align}
	which need to be equal to $-f_{1}+g_{1}$ so we need (for all $e\in Ker(d_{0})$):
	\begin{equation*}
	d_{2} \Big( \big(-s_{0}d_{0}h_{0}^{1}+h_{0}^{1}-h_{1}^{1}+s_{1}d_{0}h_{0}^{1} \big) (e) \Big)=0
	\end{equation*}
	
	On the other hand, we know that:
	\begin{equation*}
	d_{2}(l)=0
	\end{equation*}%
	for all $l\in NB_{2}$; since the Moore complex is with length one so that $NB_{2}=0$.
	
	\bigskip
	
	Now we just need to show that:
	\begin{equation*}
	\big(-s_{0}d_{0}h_{0}^{1}+h_{0}^{1}-h_{1}^{1} + s_{1} d_{0}h_{0}^{1} \big) (e)\in NB_{2}=Ker( d_{0})\cap
	Ker( d_{1})
	\end{equation*}%
	for all $e\in Ker(d_{0})$.

    \bigskip
	
	In this case:
	\begin{align*}
	d _{0} \Big( \left(-s _{0}d
	_{0}h_{0}^{1}+h_{0}^{1}-h_{1}^{1}+s _{1}d _{0}h_{0}^{1}\right) (e) \Big) & = 
	 -d _{0}s _{0}d
	_{0}h_{0}^{1}(e)+d _{0}h_{0}^{1}(e)-d _{0}h_{1}^{1}(e)+d _{0}s _{1}d
	_{0}h_{0}^{1}(e) \\
	& =  -d_{0}h_{0}^{1}(e)+d_{0}h_{0}^{1}(e)-h_{0}^{0}d_{0}(e)+s _{0}d _{0}d
	_{0}h_{0}^{1}(e) \\
	& =  -h_{0}^{0}d_{0}(e)+ s _{0}d _{0}d _{1}h_{0}^{1}(e)
	\\
	& =  -h_{0}^{0}d_{0}(e)+ s _{0}d _{0}d _{1}h_{1}^{1}(e)
	\\
	& = -h_{0}^{0}d_{0}(e)+ s _{0}d _{0}d _{0}h_{1}^{1}(e)
	\\
	& =  -h_{0}^{0}d_{0}(e)+ s _{0}d _{0}h_{0}^{0}d_{0}(e) \\
	& =  0 \ \ \ \ \ \ \ (\because e\in Ker(d_{0}))
	\end{align*}
	means:
	\begin{equation}\label{ker-1}
	\left(-s_{0} d_{0}h_{0}^{1}+h_{0}^{1}-h_{1}^{1}+ s_{1}
	d_{0}h_{0}^{1}\right) (e)\in Ker( d_{0})
	\end{equation}
	
	Similarly:
	\begin{align*}
	d _{1} \Big( \left(-s _{0}d
	_{0}h_{0}^{1}+h_{0}^{1}-h_{1}^{1}+s _{1}d _{0}h_{0}^{1}\right) (e) \Big) & = 
	-d _{1}s _{0}d
	_{0}h_{0}^{1}(e)+d _{1}h_{0}^{1}(e)-d _{1}h_{1}^{1}(e)+d _{1}s _{1}d
	_{0}h_{0}^{1}(e) \\
	& = -d _{0}h_{0}^{1}(e)+d _{1}h_{0}^{1}(e)-d
	_{1}h_{1}^{1}(e)+d _{0}h_{0}^{1}(e) \\
	& = -d _{0}h_{0}^{1}(e)+d _{1}h_{1}^{1}(e)-d
	_{1}h_{1}^{1}(e)+d _{0}h_{0}^{1}(e)  \\
	& = -d _{0}h_{0}^{1}(e)+d _{0}h_{0}^{1}(e)  \\
	& = 0%
	\end{align*}
	means:
	\begin{equation}\label{ker-2}
	\left(- s_{0} d_{0}h_{0}^{1}+ h_{0}^{1}-h_{1}^{1}+ s_{1}d_{0}h_{0}^{1}\right) (e)\in Ker( d_{1})
	\end{equation}
	Following from \eqref{ker-1} and \eqref{ker-2}:
	\begin{equation*}
    \left(-s_{0} d_{0}h_{0}^{1}+h_{0}^{1}-h_{1}^{1}+ s_{1}d_{0}h_{0}^{1}\right) (e) \in NB_{2}=Ker( d_{0})\cap Ker( d_{1})
	\end{equation*}
	
    Therefore:
	\begin{align*}
	 d_{2} \Big( \left(-s_{0} d_{0}h_{0}^{1}+h_{0}^{1}-h_{1}^{1}+ s_{1}d_{0}h_{0}^{1}\right) (e) \Big)=0
	\end{align*}
	
	Finally if we continue the calculations \eqref{continue} we get:
	\begin{align*}
	H \circ \d = -f_{1} + g_{1}
	\end{align*}
	Therefore for all $e \in Ker(d_0)$:
	\begin{align*}
		g_{1} (e)= f_{1}(e)+(H\circ \d)(e).
	\end{align*}
	
	\textbf{(iii)}  Finally, we need to show that the map $H$ satisfies the required $f_{0}$-derivation conditions given in \eqref{s-der}.
	
	\medskip
	
	The first condition is:
	\begin{align*}
		H(r+r') & = (-s _{0}f_{0} + h_{0}^{0})(r + r') \\
		& = -s _{0}f_{0} (r+r') + h_{0}^{0} (r+r') \\
		& = - \Big( s _{0}f_{0} (r) + s _{0}f_{0} (r') \Big) + h_{0}^{0} (r) + h_{0}^{0} (r') \\
		& = -s _{0}f_{0} (r') -s _{0}f_{0} (r) +   h_{0}^{0} (r) + h_{0}^{0} (r') \\
		& = -s _{0}f_{0} (r') -s _{0}f_{0} (r) +   h_{0}^{0} (r) + s _{0}f_{0} (r') - s _{0}f_{0} (r') +  h_{0}^{0} (r') \\
		& = f_0(-r') \ast \Big( -s _{0}f_{0} (r) +   h_{0}^{0} (r) \Big) - s _{0}f_{0} (r') +  h_{0}^{0} (r') \\
		& = \Big( f_0(-r') \ast H(r) \Big) + H(r'),
	\end{align*}
	
	\medskip
	
	and the second one to satisfy is:
	\begin{equation*}
	H(r \ast r')=f_{0}(r) \ast H(r') + f_{0}(r') \ast^{\circ} H(r) + H(r) \ast H(r')
	\end{equation*}
	
	On the left hand side we have:
	\begin{align*}
	H(r \ast r') & = 	(-s _{0}f_{0} + h_{0}^{0})(r \ast r')  \\
    & = -s _{0}f_{0}(r \ast r') + h_{0}^{0}(r \ast r')\\
	& = -s _{0}f_{0}(r) \ast s_{0}f_{0}(r') + h_{0}^{0}(r) \ast h_{0}^{0}(r')	 \\
	\end{align*}
	while on the right hand side is:
	\begin{align*}
	f_{0}(r) \ast
	(-s _{0}f_{0} + h_{0}^{0})(r') & + f_{0}(r')\ast (-s _{0}f_{0}+h_{0}^{0})(r)+(-s
	_{0}f_{0} + h_{0}^{0})(r) \ast (-s _{0}f_{0} + h_{0}^{0})(r') \\
	& =s _{0}f_{0}(r) \ast (-s _{0}f_{0}+h_{0}^{0})(r')+s _{0}f_{0}(r') \ast (-s _{0}f_{0}+h_{0}^{0})(r) \\
	& \quad \quad +(-s _{0}f_{0}+h_{0}^{0})(r) \ast (-s _{0}f_{0}+h_{0}^{0})(r') \\
	& =s _{0}f_{0}(r) \ast h_{0}^{0}(r')-s _{0}f_{0}(r) \ast s
	_{0}f_{0}(r')+s _{0}f_{0}(r') \ast h_{0}^{0}(r)-s
	_{0}f_{0}(r') \ast s _{0}f_{0}(r) \\
	& \quad \quad +h_{0}^{0}(r) \ast h_{0}^{0}(r')-h_{0}^{0}(r) \ast s
	_{0}f_{0}(r')-s _{0}f_{0}(r) \ast h_{0}^{0}(r')+s
	_{0}f_{0}(r) \ast s _{0}f_{0}(r') \\
	&  = -s _{0}f_{0}(r) \ast s_{0}f_{0}(r') + h_{0}^{0}(r) \ast h_{0}^{0}(r').
	\end{align*}
	which completes the proof.
	
	\medskip
	
	Remark that, in the previous calculations we explicitly used the Axiom 1, simplicial identities \eqref{simpident} and the simplicial homotopy identities \eqref{homident}.
\end{pf}

\medskip

Moreover we can give the following theorem as a consequence of the previous one:

\begin{thm}\label{main2}
	The functor $X_1$ preserves the homotopy equivalence. On other words, if $\A$ and $\B$ be simplicial objects with Moore complex of length one such that $\A \simeq \B$, then: $$ X_1 (\A) \simeq X_1 (\B) .$$
\end{thm}

\begin{pf}
	It follows from Theorem \ref{main1} and the functorial properties of $X_1$.
\end{pf}

\section{Applications}

If we handle the category $\C$ as the category of groups, which is a MCI, we get the formula of the derivation given in \cite{ JFM1, GM1} as: 
\begin{align*}
	s(gh)= \big( f_0(h^{-1}) \l s(g) \Big) \, s(h).
\end{align*}

Now let us examine the homotopies in the category of crossed module morphisms in the category of associative (bare) algebras, Leibniz algebras, Lie algebras and dialgebras (diassociative algebras) which are the examples of MCI. We refer \cite{ casasleib, CL1, dedecker1966} to recall these structures. In these constructions, we use the different types of the symbol $\l$ to denote the possible actions in such categories. Additionally, all algebras will be defined over a fixed commutative ring $\kappa$.

\subsection{Associative Algebras}

\begin{defn}
	Let $f_{0}\colon R\to R^{\prime }$ be an associative algebra (or bare algebra \cite{JFM2}) homomorphism. An $%
	f_{0}$-derivation $s\colon R\to E^{\prime }$ is a $\k$-linear map satisfying, for all $a,b \in R$:
	\begin{align}\label{s-bare}
	s(ab)=f_{0}(a) \l s(b)+ s(a) \r f_{0}(b) + s(a)s(b).
	\end{align}
	Remark that, this formula is the generalization of the derivation formula, which given for commutative algebras in \cite{IJK}.
\end{defn}

\subsection{Leibniz Algebras}

\begin{defn}
	Let $f_{0}\colon R\to R^{\prime }$ be a Leibniz algebra homomorphism. An $%
	f_{0}$-derivation $s\colon R\to E^{\prime }$ is a $\k$-linear map satisfying, for all $a,b \in R$:
	\begin{align}\label{s-leibniz}
	s \big( \ll a,b \rr \big) = f_{0}(a) \l s(b) + s(a) \r f_{0}(b) + \ll s(a),s(b) \rr.
	\end{align}
\end{defn}

\subsection{Lie Algebras}

\begin{rem}
	The notion of the homotopy of crossed modules of Lie algebras is obtained by reducing from Leibniz algebras in the sense of \cite{CL1}. Therefore the $s$ derivation formula will be (for all $a,b \in R$):
		\begin{align}\label{s-lie}
		s \big( \left[ a,b \right] \big) = f_{0}(a) \l s(b) - f_{0}(b) \l s(a) + \left[ s(a),s(b) \right].
		\end{align}
\end{rem}

\subsection{Dialgebras}

\begin{defn}
	Let $f_{0}\colon R\to R^{\prime }$ be a dialgebra homomorphism. An $%
	f_{0}$-derivation $s\colon R\to E^{\prime }$ is a $\k$-linear map satisfying, for all $a,b \in R$:
	\begin{align}\label{s-dialg}
	\begin{split}
	s(a \vdash b)&=f_{0}(a) \la s(b)+s(a) \ra f_{0}(b)+s(a) \vdash s(b), \\
	s(a \dashv b)&=f_{0}(a) \lb s(b)+s(a) \rb f_{0}(b)+s(a) \dashv s(b).
	\end{split}
	\end{align}
\end{defn}

\begin{thm}
	Let $f=(f_1,f_0)$ be any crossed module morphism $\X \to X'$ of one the categories such as associative algebras, Leibniz algebras, Lie algebras and dialgebras. In the conditions of
	previous definitions, if we define $%
	g=(g_{1},g_{0})$ as:  
	\begin{align*}
	&g_{0}(r)=f_{0}(r)+(\partial ^{\prime } \circ s)(r), 
	&&g_{1}(e)=f_{1}(e)+(s \circ \partial )(e),
	\end{align*}
	(where $e\in E$ and $r\in R$). Therefore $g$ also defines a crossed module morphism $\X \to \X'$ and we get the homotopy $f\xra{(f_{0},s)}g,$
    connecting $f$ to $g$ (Definition \ref{homotopy}). 
\end{thm}

\begin{cor}
	In the condition of all homotopy definitions given in \eqref{s-bare}, \eqref{s-leibniz}, \eqref{s-lie} and \eqref{s-dialg}; one can see that the adjoint functors between the categories: 
		\begin{align*}
		\xymatrix@C=60pt@R=70pt{ &  \mathbf{XDiAs} 
			\ar@{}[dl]^(.1){}="a"^(.9){}="b" \ar^{XAs} "a";"b"
			\ar@{}[dr]^(.1){}="a"^(.9){}="b" \ar_{XLb} "a";"b"
			&   \\
			\mathbf{XAs} \ar@<-1ex>[r]_{XLie_1}
			\ar@<2ex>@{}[ru]^(.1){}="a"^(.9){}="b" \ar^{\subset} "a";"b"
			& \mathbf{XLie} \ar@<-1ex>[l]_{XU} 
			\ar@<-1ex>[r]_{\subset}
			&  \mathbf{XLb}
			\ar@<-1ex>[l]_{XLie_2}
			\ar@<-2ex>@{}[lu]^(.1){}="a"^(.9){}="b" \ar_{U_d} "a";"b"
        }
		\end{align*}
	not only preserving the crossed module structure, also preserving the homotopy relations for crossed module morphisms in the sense of \cite{CL1}.
\end{cor}

\bibliographystyle{plain}
\bibliography{references}

\end{document}